\documentclass[12pt]{article}
 \usepackage[english]{babel}
\usepackage[T2A]{fontenc}
\usepackage{amssymb}
\usepackage{amsmath}
\usepackage{amscd}
\usepackage{setspace}
\begin{document}
	\newtheorem{mydef}{Definition}[section]
	\newtheorem{theorem}{Theorem}[section]
	\newtheorem{seqv}{Corollary}[section]
	\newtheorem{lemm}{Lemma}[section]
	\newtheorem{examp}{Example}[section]
	\newtheorem{abst}{Abstract}[section]
	\newtheorem{remark}{Remark}[section]
	\newenvironment{proof}
	{\par\noindent{\bf Proof.\enspace}}
	{\hfill$\square$}
	\author{Z.R.~Gabidullina}
	\title{A Fully Adaptive Steepest Descent Method
		\thanks{Kazan Federal University,\, e-mail: zgabid@mail.ru,\, zulfiya.gabidullina@kpfu.ru }}
	\date{}
	\maketitle
	%
%	\title{A Fully Adaptive Steepest Descent Method \thanks{}}
%	\author{Gabidullina Z.R. \inst{1}\orcidID{0000-0003-3573-3979}}
	%
%	\authorrunning{Gabidullina Z.R.}
%	\institute{Kazan Federal University,  Kazan, Russia
%		\email{zulfiya.gabidullina@kpfu.ru}\\
%		\url{http://kpfu.ru/main?p\_id=10082}}
	%
	\maketitle              % typeset the header of the contribution
	\abstract
		For solving  pseudo-convex  global  optimization problems, we present a novel  fully adaptive steepest descent method (or ASDM) without any hard-to-estimate  parameters.  For the step-size regulation in
		an $\varepsilon$-normalized direction, we use the deterministic rules, which were proposed in J. Optim. Theory Appl. (2019,\, DOI: 10.1007/S10957-019-01585-W).
		We obtained the optimistic 
		convergence  estimates for the generated by ASDM sequence of  iteration points. Namely, the sequence of function values of iterates has the
		advantage of the strict monotonic behaviour and globally converges  to the objective function optimum with the sublinear  rate.
		This rate of convergence is now known to be the best for the steepest descent method in the non-convex objectives context. 
		Preliminary computational tests confirm the efficiency of the proposed method and low computational costs for its realization.\\

		{\bf keywords:} {pseudoconvex function, steepest descent, normalization of descent direction, adaptive step-size, rate of convergence\\
{\bf MSC classes:} 90C30,\,  65K05

	%
	%
	%
	%%%%\section{First Section}
	\section{Introduction}\label{Gabid_sec1}
	As is broadly known, the development of the original variant of the
	steepest descent met\-hod (or, briefly, SDM) was
	pioneered by Cauchy (1847) for solving systems of homogeneous 
	equations. 
	For solving the mathe\-ma\-ti\-cal programs in the other  settings,
	a~wide spectrum of various kinds of SDM was investigated
	by re\-sear\-chers over the  many years.  A~very useful systematic
	survey of the existing literature related to the different variants
	of SDM can be found, for example, in \cite{Gabid_Burachik,Gabid_Asmundis,Gabid_Bento_etc} (see also
	the references therein). For an up-to-date survey of the topic, it is
	appropriate to refer to \cite{Gabid_Asmundis}.\,
	
	To avoid some conflicts and ambiguities that can be arisen in what the term "steepest descent method" means in the literature on optimization framework, we will try here to clarify some terminology. A question can now be addressed. What conditions on an optimization method ensure that the name SDM is  properly used for its characterizing. The reasons why we have felt the need for such an explanation are twofold. In a wide range of optimization topics, some gradient method is usually called a steepest descent method in the case of using the uniform descent indicated by opposite to the gradient direction.  Although tradition merely associates SDM not only with specific descent directions but with selecting the step-size by exact (or so-called perfect) line search. At the same time, the wide spectrum of papers on global optimization  utilize the term "steepest descent method" or modifications of SDM regardless of what strategies is used for the step-size selection. In this paper, we  apply in ASDM the rules of regulating the step length which are different from the exact (computationally expensive) line search.
	
	Rather than describing all the various versions of SDM, that
	researchers have been constructed over the years for the achievement of
	the best results, we will occupy our attention only in some
	relevant work.
	
	In \cite{Gabid_	Goldstein},\cite{Gabid_	Goldstein_Price},\,  there were studied the two versions of SDM for functions being twice differentiable  on the Euclidean space. The second of them provides global convergence at a rate which is eventually superlinear and possibly quadratic.
	
	In \cite{Gabid_Vrahatis_etc}, the development, convergence theory and numerical testing some versions of steepest descent algorithm with adaptive step-size was presented. All of the algorithms are computationally efficient. Based on estimates of the Lipschitz constant, there was proved the convergence of the different variants of SDM to a minimizer or to a stationary point of objective function. The algorithms have been tested on real-life artificial neural network applications and the results have been very satisfactory.
	
	To minimize a continuously differentiable quasiconvex function, 
	SDM with Armijo's step-sizes was proposed in \cite{Gabid_KIWIEL_MURTY}. This method generates a sequence of iterates  globally converging to a point at which the gradient of the objective function is equal to zero vector.
	
	In \cite{Gabid_Barzilai_Borwein}, there were derived two-point step sizes for SDM by approximating the secant equation. At the cost of storage of an extra iterate and gradient, these algorithms achieve better performance and cheaper computation than the classical SDM. By Barzilai and Borwein, there was proposed the non-monotonic variant of SDM  which is superlinearly convergent for convex quadratic setting in two-dimensional space, and has quite well behaviour for the case of high-dimensional tasks. 
	
	In \cite{Gabid_Fletcher}, there is used an idea that a limited memory approach
	might be fashioned by using a limited number of eigenvalue estimates. An improvement of characteristics of  SDM  has been achieved
	by the introduction of the Barzilai-Borwein choice of step length, and some other related ideas. There is suggested a method which is able to take advantage of the availability of a few additional `long' vectors of storage to achieve a significant improvement in performance, both for quadratic and non-quadratic objective functions. The sequence of iterates converges to the point for which the gradient of objective function equals the zero vector.
	
	A step-size formulae, which provides for SDM fast convergence and the monotone property, was presented in  \cite{Gabid_Yuan}. An algorithm with the new step-size in even iterations and exact line search in odd iterations is proposed. Numerical results obtained by the new method confirm  that the the exact solution may be found  within three iterations for the case of two-dimensional problems. For small-scale problems, the new method is very efficient.  A modified version of the new method is also presented, where the new technique for choosing the step length is utilized after every two iterations with exact line searches. The modified algorithm is comparable to the Barzilai-Borwein method for large-scale tasks and better for small-scale tasks.
	
	In \cite{Gabid_SAHU}, there is investigated a generalized hybrid steepest descent method and its convergence theory for solving monotone variational inequality over the fixed point set of a mapping which is not necessarily Lipschitz
	continuous. This method is used for solving the convex minimization problem for a smooth convex function whose gradient is not necessarily Lipschitzian. There is proved that the sequence of iterates converges strongly to  a minimizer \,$x^{*}$.\,
	
	Full convergence of the steepest descent method with inexact line searches was proved in \cite{Gabid_Burachik}. There were  considered two of such procedures and proved, for a convex objective function, convergence of the whole sequence to a minimizer without any level set boundedness assumption and, for one of them, without any Lipschitz condition.
	
	In \cite{Gabid_Asmundis}, there was demonstrated how taking into account the spectral properties of the Hessian matrix, for convex quadratic problems, one can provide the improvement of the practical behaviour of SDM. This allows them
	to obtain computational results comparable with those of the Barzilai and Borwein algorithm, with the further advantage of monotonic behaviour.
	
	In \cite{Gabid_Nocedal},\, there are  presented results describing the properties of the gradient norm for the SDM applied to quadratic objective functions. There are also made some general observations that apply to nonlinear problems, relating the gradient norm, the objective function value, and the path generated by the iterates.
	
	By the method from \cite{Gabid_Bento_etc}, there is guaranteed the well definedness of the  generated sequence. Under mild assumptions on the multicriteria function, there was justified that each accumulation point (if they exist) satisfies first-order necessary conditions for Pareto optimality. Under assumptions of quasi-convexity of the multicriteria function and non-negativity of
	the Riemannian manifold  curvature, full convergence of the sequence to a Pareto critical was proved.

	The main contributions in this paper are as follows. We propose
	a~novel fully adaptive steepest descent method with the step length regulation for solving  pseudo-convex unconstrained optimization tasks. This relaxation al\-go\-rithm
	allows one to generate the sequence of iterates \,$\{x_{k}\}$,\,
	\,$k = 0,\, 1, \ldots$\, such that the sequence of its function
	values \,$\{f(x_{k})\}$,\, \,$k = 0,\, 1, \ldots$\, has the
	advantage of the strict monotonic behaviour and converges
	globally to the objective function optimum with the
	following rate \,$O (1 /k ).$\,  This rate is traditionally called the
	sublinear one. To the best of our knowledge, a convergence rate
	of \,$O (1 /k )$\, is now known to be the best for SDM  in a~non-convex objectives context.
	
	The sublinear convergence rate takes place
	under the following conditions relating the original problem: 1) an objective function \,$f(x)$\, is  pseudo-convex on some convex set \,$D \subseteq \mathbb{R}^{n}$ (the set \,$D$\, may coincide, for instance, with the Lebesgue set (corresponding to a starting point of the iterates sequence) of the objective function or with the whole Euclidean space and etc),\, 2) the function \,$f(x)$\, is required
	to be satisfied to  so-called Condition~A\, introduced in
	\cite{Gabid_gab1}. We note that this condition will be defined explicitly
	below in Section~2 (see De\-fi\-ni\-tion~\ref{Gabid_condA}).
	
	Here, we underline that, for the execution of ASDM, there is no need to make use any priori information regarding the auxiliary constant defined in Condition~A.
	With respect to this fact, the presented version of ASDM compares favorably with the
	other variants of SDM considered above.
	The fulfillment of Condition~A allows us further to adaptively regulate the step-size
	without appling  any complicated line search techniques. Indeed, there is no need to
	use the one-dimensional exact minimization of the objective function in the selected descent direction.
	We apply the two deterministic rules of the step length adjustment. These strategies guarantee to determine  the step-size by utilizing the finite procedures of diminishing an initial value of a certain parameter.
	The latter is a user-selected parameter which is diminished until
	a moment when the condition applied for the step length regulation
	becomes fulfilled. We notice that the step-size computing
	strategies provide the strict relaxation of the objective function
	at each iteration.
	Note that the concept of the objective function relaxation
	allows to interpret of the relaxation properties of
	minimization methods. Due to this interpretation,  we can to
	evaluate how the objective function value is decreased at
	each iteration. The event of this value being diminished on the
	positive magnitude corresponds to the property of the strict
	relaxation which is described by the following inequality: \,$ f(x_{k}) >
	f(x_{k+1})$,\, $k = 0,\, 1,\, \ldots$\,
	
	The fully adaptive character of the presented variant of
	SDM is established na\-me\-ly
	by combining the simultaneous
	control of adapting 
	an~$\varepsilon-$normalization parameter
	of the descent direction as well as the step-size
	re\-gu\-la\-tion in tandem.
	We establish the finiteness of all the procedures for the
	step-size regulation as well as the adaptation of the
	$\varepsilon-$normalization parameter.
	
	Due to all the mentioned properties of ASDM, it seems that the
	results of the paper may be potentially useful in various applied domains covering 
	their theoretical as well as practical aspects (in
	pseudo-convex prog\-ram\-ming, variational inequality problem solving, and many others
	(in particular, data clas\-si\-fi\-ca\-tion 
	techniques and neural networks simulation)).
	
	The rest of this paper is organized as follows. In
	Section~2 we present some preliminaries which are necessary for our convergence rate
	analysis of ASDM. In Sec\-tion~3, there is formulated ASDM and justified
	its convergence rate. In Section~4, there are
	drawn some conclusions. 
	
	\section{Definitions and Preliminaries}\label{Gabid_sec2}
	
	In this paper we aim to to explore the following problem:
	\begin{equation} \label{Gabid_eq1} \min  \limits_{x \in \mathbb{R}^{n}}\, f(x),
	\end{equation}
	where \,$f(x):\mathbb{R}^{n} \rightarrow \mathbb{R}^{1}$\, is a continuously differentiable pseudo-convex
	function sa\-tis\-fying  the so-called Condition \,$A$\,
	(introduced in \cite{Gabid_gab1}) on a convex set
	\,$D \subseteq \mathbb{R}^{n}$.\,
	To solve this problem, we propose a new efficient algorithm, which
	has the estimates of the rate of its convergence and allows one
	to adaptively re\-gulate both the parameter of
	an~$\varepsilon-$nor\-ma\-lization of a descent direction and a~step length.
	
	We begin with some notations: 
	\,$$\nabla f(x) = (\dfrac{\partial f(x)}{\partial
		\xi_{1}},\, \dfrac{\partial f(x)}{\partial
		\xi_{2}},\,\ldots,\,\dfrac{\partial f(x)}{\partial
		\xi_{n}})$$
	is the gradient of the
	function \,$f(x)$\, at the point \,$x = (\xi_{1},\,\xi_{2},\, \ldots, \xi{n})$,\, \,$x_{0}$\, stands for
	a starting point of the iterative con\-se\-qu\-ence \,$\{x_{k}\},\, k \in \mathbb{N}$\, generated by
	minimizing the objective function.
	
	Let \,$\|\cdot\|$\, stand for
	the Euc\-li\-de\-an norm  of a vector in $\mathbb{R}^{n}$,\,  \,$\langle \cdot, \cdot \rangle$\, stand for the usual inner product,\,   \,$f^{*}:= \min \limits_{x \in
		\mathbb{R}^{n}}\, f(x)$,\, \,$X^{*} := \{ x \in \mathbb{R}^{n}: f(x) = f^{*}\}$,\,
	\,$\mathbb{N} = \{0,\,1,\,\ldots\}$,\, ${\bf 0}$ be a zero vector of \,$\mathbb{R}^{n}$,\,  and \,$p_{k}^{*}$\,
	correspond to a~pro\-jec\-tion of the iterative point \,$x_{k}$\,
	onto the set \,$X^{*}$,\, $k \in \mathbb{N}$.\,  In the literature on optimization, \,$p_{k}^{*},\, k \in \mathbb{N}$\, are sometimes called accumulation points.
	
	To the extent of our knowledge,  the class of
	smooth pseudo-convex functions was pioneered
	by Mangasarian in \cite{Gabid_mang}. The
	above-mentioned class represents a~generalization of the family of all continuously differentiable  convex
	functions.
	
	\begin{mydef} \rm{(pseudo-convexity)
			A  function \,$f(x)$,\, which is given and continuously differenti\-ab\-le on
			an~open and convex set \,$G$\, from \,$\mathbb{R}^{n}$,\, is called
			pseudo-convex, if there is fulfilled the
			following implication:
			$$\langle \nabla f(x), y - x \rangle \geq 0 \Rightarrow f(x) \leq f(y),\,\, \forall  x,\, y \in G,$$
			or equivalently, $$f(y) < f(x) \Rightarrow \langle \nabla f(x), y - x
			\rangle < 0,\,\, \forall  x,\, y \in G.$$}
	\end{mydef}
	In the case of pseudo-convex functions,  the necessary and
	sufficient con\-di\-tions of optimality are established in the
	following theorem.
	\begin{theorem} \label{Gabid_t2.1} (basic first-order conditions for
		optimality) \rm{(\cite{Gabid_mang}, p.282)}  \it{For the point
			\,$x^{*}\in G$\, to furnish the minimum of \,$f(x)$\, over
			\,$G$,\, it is necessary and sufficient for all \,$x \in G$\, to
			hold}
		$$\langle \nabla f(x^{*}), x - x^{*} \rangle \geq 0.$$
	\end{theorem}
	\begin{mydef} \label{Gabid_condA}  \rm{(Con\-di\-tion~$A$)
			We say that a continuous function \,$f(x)$\, sa\-tis\-fi\-es
			Con\-di\-tion~$A$ on the convex set \,$D \subseteq
			\mathbb{R}^{n}$\, if there exist a nonnegative symmetric function
			\,$\tau(x,y)$\, and \,$\mu > 0$\,  such that}
		\begin{multline}\nonumber
		f(\alpha x + (1 - \alpha)y ) \geq \alpha f(x) + (1 - \alpha)f(y) -
		\alpha(1 -\alpha)\mu \tau(x,y),\\ \forall x,\, y \in D,\,\, \alpha
		\in [0, 1].
		\end{multline}
	\end{mydef}
	For \,$x,\, y \in D\subseteq \mathbb{R}^{n}$,\, we say that some function
	\,$\tau(x,y)$\, is symmetric if \,$\tau(x,y) =
	\tau(y,x)$,\, $\tau(x,x) = 0$.\, Condition~$A$ characterizes
	a~sufficiently wide class of functions \,$A(\mu,\tau(x,y))$\,. It
	was demonstrated in
	\cite{Gabid_gab1},\cite{Gabid_gab2,Gabid_gab3} that the functions 
	class \,$A(\mu,\|x-y\|^{2})$,\, in
	particular,
	is broader than \,$C^{1,1}(D)$\, - the  commonly known class of functions
	whose gradient vectors have Lipschitzian property on the convex set
	\,$D \subseteq \mathbb{R}^{n}$.\, By the way, we notice that namely
	the Lipschitz condition for gradients of functions being minimized
	has  been determined as the necessary assumption in justifying 
	the theoretical estimates of the convergence rate for many
	modern smooth optimization methods. In
	\cite{Gabid_gab3}, there were presented some examples  of
	functions that satisfy Condi\-ti\-on~A. For functions from
	\,$A(\mu,\tau(x,y))$,\, there also were explored  their main
	properties and criteria for membership in the studied class.
	Furthermore, for a smooth function satisfying Condition A on a convex
	set D,  there was proved in \cite{Gabid_gab3}  the following remarkable differential inequality:
	\begin{equation}\label{Gabid_e1.10}
	f(x) - f(y) \geq \langle \nabla f(x), x - y \rangle - \mu \tau(x,y).\,
	\end{equation}
	\begin{theorem} \label{Gabid_t1.5_S1.1} (relation between two classes of functions) \rm{(\cite{Gabid_gab11}, p.1082)}
		\it{If \,$D$\, is convex subset of \,$\mathbb{R}^{n}$,\, \,$f(x) \in
			C^{1,1}(D)$,\, then \,$f(x)$\, satisfies  Condition A  on
			\,$D$\, with coefficient \,$\mu = L/2$\, and
			function \,$\tau(x,y) = \|x - y\|^{2}$,\, where \,$L$\, is
			a Lipschitz constant for the gradient of  \,$f(x)$.}
	\end{theorem}
	
	\begin{mydef}  \rm{($\varepsilon -$normalized descent direction)\enspace
			For functions from the class
			\,$A(\mu,\|x-y\|^{v})$,\, \,$v \geq 2$,\,  a vector \,$s \neq {\bf 0}$\,
			is called an $\varepsilon-$normalized descent direction
			($\varepsilon > 0$) of the function $f$ at the point $x \in D \subseteq
			\mathbb{R}^{n}$ if
			the following inequality is fulfilled:}
		$$\langle \nabla f(x),s \rangle + \varepsilon \|s\|^{v} \leq 0.$$
	\end{mydef}
	\begin{lemm}\label{Gabid_lem3.1}(\,$\varepsilon -$normalization)
		If some descent direction  \,$s$\, is not
		\,$\varepsilon -$nor\-ma\-li\-zed, then the vector constructed in
		such a way that \,$\bar{s} = \genfrac{}{}{}{0}{t
			s}{\varepsilon\|s\|^{v}}$\, is \,$\varepsilon-$nor\-ma\-li\-zed
		under the condition $0 < t \leq |\langle \nabla f(x),s \rangle|$.\,
	\end{lemm}
	\proof By construction, we obviously obtain the following relation:
	\begin{multline}\nonumber
	\langle \nabla f(x),\bar{s} \rangle +
	\varepsilon \|\bar{s}\|^{v} =
	\genfrac{}{}{}{0}{t}{\varepsilon\|s\|^{v}} \langle \nabla f(x),s \rangle
	+ \genfrac{}{}{}{0}{t^{v}}{\varepsilon^{v - 1}\|s\|^{v(v - 1)}}
	=\\
	=  \genfrac{}{}{}{0}{t}{\varepsilon\|s\|^{v}}\left[\langle \nabla f(x),s
	\rangle + t \left[\genfrac{}{}{}{0}{t}{\varepsilon\|s\|^{v}}
	\right]^{v - 2} \right] \leq 0,
	\end{multline}
	because \,$\left[\genfrac{}{}{}{0}{t}{\varepsilon\|s\|^{v}}
	\right]^{v - 2} \leq 1.$\, $\square$
	\smallskip
	
	Under the condition $v = 2$, if 
	for some fixed point \,$x \in
	\mathbb{R}^{n}$\, the vectors \,$z - x$\, are
	$\varepsilon -$normalized descent directions at the point
	$x$,\, then it is easy to observe that all the points
	\,$z \in \mathbb{R}^{n}$\, belong to the $n-$dimensional ball  of radius \,$R =
	\|\nabla f(x)\|/2\varepsilon$\, with center at the point $u = x -
	\nabla f(x)/2\varepsilon$.\,
	
	Set
	\[\zeta = \begin{cases}
	(\varepsilon\cdot\mu^{-1})^{1/(v-1)},\, & \text{if \,$\varepsilon < \mu,$}\\
	1,\, & \text{if \,$\varepsilon \geq \mu.$}\,
	\end{cases}\]
	We
	further study very useful properties of \,$\varepsilon-$normalized descent di\-rec\-tions. These properties  provide a
	strict relaxation of the objective function in gradient methods.
	\begin{lemm} \label{Gabid_lem3.4} \,(main properties of \,$\varepsilon-$nor\-ma\-li\-zed descent directions) (\cite{Gabid_gab11}, p.1083)\, 
		Let \,$s$\, be some \,$\varepsilon-$nor\-ma\-li\-zed descent direction  for
		the function \,$f$\,
		at the point \,$x$\, where \,$v \geq 2$,\, \,$f(x) \in A(\mu,\|x-y\|^{v})$,\, then for all \,$\beta \in \left] 0,\,1\right[ $\,  there
		exists a constant \,$\hat{\lambda} = \hat{\lambda}(\beta) > 0$\,
		\,$(\hat{\lambda} = (1 - \beta)^{1/(v-1)}\zeta)$\, such that
		for all \,$\lambda \in \left] 0,\,\hat{\lambda} \right] $\, it
		holds
		\begin{equation}\label{Gabid_eqq1}
		f(x) - f(x+\lambda s)
		\geq -\lambda \beta \cdot \langle \nabla f(x), s \rangle,
		\end{equation}
		\begin{equation}\label{Gabid_eqq2}
		f(x) - f(x+\lambda s) \geq \lambda \beta
		\cdot \varepsilon \|s\|^{v}.
		\end{equation}
	\end{lemm}
	For establishing the convergence of the adaptive algorithm,
	which will be proposed below, there are
	needed the inequalities (\ref{Gabid_eqq1})--(\ref{Gabid_eqq2}).\, In particular, they imply a strict relaxation of the objective function. Namely, it holds $$f(x+\lambda s) < f(x),\,  \,\forall\, \lambda \in (0,\,  (1 - \beta)^{1/v-1}\zeta],\, \beta
	\in (0,1).$$
	Based on Lemma \ref{Gabid_lem3.4},\,
	there can be described the rules of computing the step-size satisfying
	(\ref{Gabid_eqq1})--(\ref{Gabid_eqq2}).\,
	Suppose that $s$ is some
	$\varepsilon-$normalized direction of descent for $f$ at the
	point $x$. Additionally, let the following conditions be fulfilled:
	\, $\beta \in \left] 0,\,1\right[ ,$\, \,$\eta = (1 -\beta)^{1/v-1},$\, \,$\hat{i} = 1,$\,
	\,$J(\hat{i}) = \{\hat{i}, \hat{i} + 1,
	\hat{i} + 2, \ldots\}.$\, Next we need to find \,$i^{*}$\, - the least
	index \,$i \in J(\hat{i})$\, for which there holds the following
	inequality:
	\begin{equation}\label{Gabid_e1}
	f(x) - f(x + \eta^{i}s)
	\geq -\eta^{i} \beta \cdot \langle \nabla f(x), s \rangle,
	\end{equation}
	or the more weak inequality:
	\begin{equation}\label{Gabid_e2}
	f(x) - f(x + \eta^{i} s) \geq \eta^{i} \beta
	\cdot \varepsilon \|s\|^{v}.
	\end{equation}
	We further set $\lambda = \eta^{i^{*}}.$  In what follows, we have in view that
	there is utilized Rule~1 (or Rule~2) when  we follow
	the first (or the second) of the above-mentioned strategies for calculating 
	the step-size.
	The step length 
	determined in accordance with these strategies satisfies
	(\ref{Gabid_eqq1}) or (\ref{Gabid_eqq2}), respectively.
	
	Further, there should be attentionally  explored the case when
	\,$s$\, is the \,$\varepsilon -$normalized  direction of descent, but it
	is not \,$\mu -$normalized (this case is possible only for
	\,$\varepsilon < \mu$).\, Under the assumption that \,$0 < \varepsilon <
	\mu$,\, for the event of choosing \,$\lambda$\, according to
	%(\ref{Gabid_e1}) or (\ref{Gabid_e2})
	Ru\-le~1 or Ru\-le~2, we demonstrate that the step length 
	is bounded from below. This evidently yields that the described
	procedures of diminishing the step-size are finite.
	\begin{lemm}\label{Gabid_lem3.6} (exact lower estimate of the
		step length)(\cite{Gabid_gab11}, p.1084)\,
		If\\
		(b)\enspace $f(x) \in A(\mu,\|x-y\|^{v})$,\, $v \geq 2$,\\
		(c)\enspace $0 < \varepsilon < \mu,\, \beta \in \left] 0,\,1\right[,$\\
		(d)\enspace $s -$\, is an $\varepsilon -$normalized descent direction of
		the function $f$
		at the point $x$,\, but it is not $\mu -$normalized,\\
		(e)\enspace $i^{*}$ is the smallest index $i = 1, 2, \ldots,$ for which
		there is fulfilled the condition of Rule~1 or Rule~2,
		%(\ref{Gabid_e1})  \,$\left( \mbox{or}\, \,(\ref{Gabid_e2})\right)$,\,
		\,$\lambda = \eta^{i^{*}}$; \\
		Then the following estimate holds: \,$$\lambda > \left(\varepsilon
		\mu^{-1}\cdot(1 - \beta)^{2}\right)^{1/(v - 1)}
		> 0.$$
	\end{lemm}
	\begin{remark} \label{Gabid_note1} \rm{(exact lower estimate of the
			constant \,$\mu$)
			Due to Lemma \ref{Gabid_lem3.6}, 
			there comes im\-me\-diately the following estimate:
			\begin{equation}\label{Gabid_eq3.10}
			\mu > \varepsilon\cdot(1 - \beta)^{2}\lambda^{1-v}.
			\end{equation}
			Later, the estimate (\ref{Gabid_eq3.10})\,
			will be utilized in ASDM for adaptive regulation of the parameter for
			the $\varepsilon -$normalization of the descent direction.}
	\end{remark}
	\section{Adaptive  Steepest  Descent Algorithm and its Convergence}\label{Gabid_sec3}
	This section is aimed at providing the principles of
	selecting the $\varepsilon
	-$nor\-malization pa\-ra\-me\-ter  for the descent
	direction. We note that these principles are universal and may be utilized for developing the various adaptive algorithms with normalized descent directions (see, for instance, \cite{Gabid_gab11}).
	
	The method con\-ver\-gence
	for the fixed parameter
	\,$\varepsilon$\, (in the case of an arbitrary ratio
	of the parameter \,$\varepsilon$\,
	and the value of \,$\mu$\, in Condition A)
	follows from the convergence
	of the adaptive variant of SDM. We notice
	that generally saying  the constant \,$\mu$\,
	is unknown beforehand.
	In practice, the selection of the \,$\varepsilon$\,
	values  close to the \,$\mu$\, value is therefore decisive
	for the al\-go\-rithm con\-ver\-gence.
	If one chooses the too small parameter \,$\varepsilon$,\,
	then, according to Rule~1 and Rule~2,
	this may imply 
	significant diminishing the step-size.
	In the case of selecting the unjustifiably large value  of \,$\varepsilon$,\,
	the convergence of the adaptive algorithm may be slowed down. Consequently, it is expedient
	to evaluate the parameter  \,$\varepsilon$\, in the process of executing the algorithm.  
	The inequalities
	(\ref{Gabid_e1})--(\ref{Gabid_eq3.10}) allow us to make
	an adjustment to the value \,$\varepsilon$\,  by increasing it if the
	previous choice was unsuccessful. Now, we specify further details
	of a~procedure for pointwise adapting the parameter
	\,$\varepsilon$\, during the iterative implementation of the algorithm.\\
	For the iterate $x_{k}$ of the adaptive algorithm, let
	$\varepsilon_{k} > 0$ be the value of a parameter for an $\varepsilon
	-$normalization of descent direction. Let $s_{k}$ be an
	$\varepsilon -$normalized descent direction of the function $f$ at
	$x_{k}$; the iterative step-size \,$\lambda_{k}$\, is chosen
	in accordance with one of the Rules~1--2.\,
	%(\ref{Gabid_e1})--(\ref{Gabid_e2})).
	Suppose that  \,$i_{k}$\, is the least
	index $i \in J(\hat{i})$, for which there is fulfilled  the condition  of
	selecting the iteration step-size for \,$x = x_{k}$,\,
	\,$\varepsilon = \varepsilon_{k}$,\, \,$s = s_{k}$.\, According to
	Lem\-ma~\ref{Gabid_lem3.4},\, if \,$i_{k} = \hat{i}$,\, then for
	implementing the next - $(k + 1)-$th iteration it is expedient to
	remain unchanged the value of the normalization parameter, i.e. to put
	\,$\varepsilon_{k + 1 } = \varepsilon_{k}$.\, During the process
	of dropping the step-size, let the verified condition
	(\ref{Gabid_e1}) (or condition (\ref{Gabid_e2})) be fulfilled for
	\,$i_{k} > \hat{i}$.\, Then, according to
	(\ref{Gabid_eq3.10}),\, there should be increased the current
	value of the parameter for the $\varepsilon -$normalization of the
	descent direction, for example, as follows:
	\,$\varepsilon_{k + 1 } = \varepsilon_{k} \cdot \zeta_{k}$.\,
	Regardless
	of what rule is chosen for computing
	the step-size, here we have
	\begin{equation}\label{Gabid_zeta_rule}
	\zeta_{k} =
	(1 - \beta)^{1-i_{k}}.
	\end{equation}
	The fact that
	\,$\mu$\, is finite implies that after the finite number of
	increases, the value of the parameter \,$\varepsilon$\, can exceed \,$\mu$\,
	and cease to vary. Suppose that \,$j > 0$\, is the number of iteration, on which it is fulfilled \,$\varepsilon_{j} \geq \mu$.\, One then has \,$\varepsilon_{k}
	\geq \varepsilon_{j} \geq \mu$,\, \,$\forall k \geq j$.\, In this
	event, after some iteration  \,($j \geq 0$\, the adaptive algorithm
	starts to be implemented with the fixed constant for the \,$\varepsilon
	-$normalization of a descent direction. We emphasize that
	beginning with the \,$j -$th\, iteration, the step-size becomes
	unchanged: \,$\lambda_{k} = \eta$,\, \,$\forall k \geq j$.\, At
	that time, for computing the iterative step length, one needs only
	one calculation of the objective function value at the point
	\,$x_{k} + \eta s_{k}$\, (for verifying the fulfillment of the
	condition applied for selecting the
	step-size).\\
	
	\smallskip
	{\bf Algorithm}. \smallskip\\
	Step 0. Initialization. Select \,$x_{0} \in \mathbb{R}^{n}$,\, \,$\beta \in
	\left] 0,\,1\right[  ,$\, \,$\varepsilon_{0}
	> 0.$\,  Set the iteration counter \,$k$\, to \,$0$.\,\\
	Step 1. For the objective function \,$f(x)$,\, calculate the gradient vector at \,$x_{k}$.\,
	Verify the optimality criterion:\\ if $\nabla f(x_{k}) = {\bf 0},$
	then terminate the al\-go\-rithm implementation   (since $x_{k}$
	is a solution of the problem (\ref{Gabid_eq1})). Otherwise,  set \,$p_{k} = - \nabla f(x_{k}),$\,
	\[s_{k} =  \begin{cases}
	\,p_{k},\, &\text{if $\langle \nabla f(x_{k}), p_{k} \rangle
		+ \varepsilon_{k}\|p_{k}\|^{v} \leq 0$,}\\
	\genfrac{}{}{}{0}{t_{k}p_{k}}{\varepsilon_{k}\|p_{k}\|^{v}}
	= \genfrac{}{}{}{0}{p_{k}}{\varepsilon_{k}\|p_{k}\|^{v-2}} ,\, &\text{othewise.}\,
	\end{cases}\]
	
	Here $t_{k} = |\langle \nabla f(x_{k}), p_{k} \rangle| = \|\nabla f(x_{k})\|^{2} .$\\
	Step 2. Let $i_{k}$ be the least index $i \in J(\hat{i})$ for
	which there holds the condition from  Rule~1 or Rule~2\, when $x = x_{k}$, $s = s_{k}$, $\varepsilon =
	\varepsilon_{k}$. Set $\lambda_{k} = \eta^{i_{k}}$.\\
	Step 3. Compute the next iterate $x_{k+1} = x_{k} +
	\lambda_{k}s_{k}.$ \\
	Step 4. Update $\varepsilon_{k + 1} = \zeta_{k}\varepsilon_{k}.$
	Set $k = k + 1$  and go to Step 1. \smallskip \smallskip
	
	Clearly, there should be applied the same rule for choosing the step-size  at
	each iteration point.
	
	\begin{remark} \label{Gabid_note2} \rm{(verification of the descent direction)}
		Let \,$\bar{s}_{k} = \|p_{k}\|.$\, The vector \,$s_{k}$\,
		generated by the algorithm  is
		the	\,$\varepsilon-$normalized descent direction. This is obvious when \,$s_{k} = \bar{s}_{k}$.\,
		In the case of \,$s_{k} = \genfrac{}{}{}{0}{t_{k}\bar{s}_{k}}
		{\varepsilon_{k}\|\bar{s}_{k}\|^{v}}$,\, Lemma \ref{Gabid_lem3.1}
		confirms that \,$s_{k}$\, is also \,$\varepsilon-$ normalized.
		In addition, it clear\-ly holds
		\,$\|s_{k}\| \leq \|\bar{s}_{k}\|.$\,
		Really,
		\[\|s_{k}\|=\begin{cases}
		\|\bar{s}_{k}\|,&\text{if $\langle \nabla f(x_{k}),\bar{s}_{k}\rangle
			+\varepsilon_{k}\|\bar{s}_{k}\|^{v} \leq 0$,}\\
		\genfrac{}{}{}{0}{t_{k}}{\varepsilon_{k}\|\bar{s}_{k}\|^{v-1}}<\|\bar{s}_{k}\|,
		&\text{otherwise,}
		\end{cases}\]
		because
		\,$\genfrac{}{}{}{0}{t_{k}}{\varepsilon_{k}\|\bar{s}_{k}\|^{v}} =
		\genfrac{}{}{}{0}{ -\langle \nabla f(x_{k}), \bar{s}_{k}
			\rangle}{\varepsilon_{k}\|\bar{s}_{k}\|^{v}} < 1.$\,
	\end{remark}
	Suppose that $\{x_{k}\}$ is the sequence generated by ASDM.
	To explore the rate convergence of ASDM for the pseudo-convex setting,
	it is necessary to establish a criterion of global optimality for a solution.
	Usually, we do not know beforehand the optimum value of the function being minimized.
	Therefore, it is crucial to work out directly the stopping criterion of ASDM for the pseudo-convex setting.
	We note that the set \,$D$\, in the formulation of the next theorem may coincide, for instance, with the Lebesgue set of the function \,$f(x)$\, at the point \,$x_{0}
	\in \mathbb{R}^{n}$\, or with the whole space.
	\begin{theorem} \label{Gabid_t11}(constructive measure of optimality for ASDM)
		Let $f(x)$ be a~con\-ti\-nuous\-ly differentiable pseudo-convex
		func\-tion on some convex set $D \subseteq \mathbb{R}^{n},$\, \,$ X^{*}\subset D.$\, Then, 
		for the fulfillment of the equality \,$f(x_{k}) = f^{*},$\, $x_{k} \in D$,\, $k = 1,\,2,\, \ldots,$\, it is sufficient to hold
		\begin{equation}\label{Gabid_eq100}
		\nabla f(x_{k}) = {\bf 0}.
		\end{equation}
	\end{theorem}
	\proof The equality (\ref{Gabid_eq100}) obviously yields that \,$\langle \nabla f(x_{k}), x - x_{k} \rangle = 0, \forall x \in D.$\, By  definition of a~pseudo-convex function, we then have \,$f(x_{k}) \leq f(x),$\, $\forall x \in D.$ Since \,$ X^{*}\subset D,$\, this is what we want to prove. \enspace $\square$

	Further, we assume that \,$x_{k} \notin X^{*},$\, \,$\forall k
	= 0, 1, \ldots$\, With the purpose of investigating  the
	con\-ver\-gence of optimization methods in the pseudo-convex
	setting, there is usually described in the literature on
	optimization an auxiliary numeric sequence $\{\theta_{k}\}$ as
	follows:
	\begin{equation} \label{Gabid_eq22}
	\theta_{k} > 0,\, 0 < \theta_{k}\cdot (f(x_{k}) - f(x^{*})) \leq
	\langle \nabla f(x_{k}), x_{k} - x^{*} \rangle,\, x^{*} \in X^{*}, k \in
	\mathbb{N}.
	\end{equation}
	By the definition of pseudo-convex functions, there must exist 
	such values $\theta_{k}$. For instance, in the case of a continuously differentiable convex function, it holds $\theta_{k} = 1, k =0, 1, \ldots$\, The
	estimates of elements
	of the sequence \,$\{\theta_{k}\}$\, were studied, for instance, in
	\cite{Gabid_gab1,Gabid_gab3}.
	
	Before formulating   the theorem  on the  convergence  of  the sequence  \,$\{x_{k}\}\,$ constructed by ASDM
	to a solution of problem~(\ref{Gabid_eq1}),\, there should be reminded the
	following well-known fact related to con\-ver\-gence of some special numeric
	sequence.
	\begin{lemm} \label{Gabid_lem2.5}(sublinear rate of convergence for numeric sequences)
		\rm{(\cite{Gabid_Vas},\, p.102)}
		\it{If a~nu\-meric sequence $\{a_{k}\}$ is such that
			$$a_{k} \geq 0,\, a_{k} - a_{k+1} \geq q \cdot a_{k}^{2},\, k = 1,\,2,\, \ldots,$$
			where $q$ is some positive constant, then  the following estimate
			holds:
			$$a_{k} \sim O(1/k),$$
			i.e. there will be found a constant \\ \smallskip \centerline
			{$q_{1}
				> 0$\, \mbox{such that} \,$0 \leq a_{k} \leq q_{1} \cdot k^{-1},\,
				k = 1, 2, \ldots$}}
	\end{lemm}
	The next auxiliary lemma is needed for the purpose of proving the convergence theorem. This lemma establishes the upper bound for the adapted values of the normalization parameter for each iteration. 
	\begin{lemm}\label{Gabid_lem3.7}(boundedness of adapted
		values of the normalization parameter) (\cite{Gabid_gab11}, p.1088)\, 
		If   \\(b1)\enspace $f(x) \in A(\mu,\|x-y\|^{v})$,\, $v \geq 2$,\\
		(c1)\enspace $s_{k}$
		is the $\varepsilon_{k}-$normalized descent
		direction,\\
		(d1)\enspace $i_{k}$ is the least index $i \in J(\hat{i})$ for which there holds
		one of the conditions (\ref{Gabid_e1}) or (\ref{Gabid_e2})  under the
		assumptions that $s = s_{k},$
		$x = x_{k},$ $\varepsilon = \varepsilon_{k},$ $\lambda_{k} =
		\eta^{i_{k}},$\\ (e1)\enspace $\{x_{k}\}$ is some iterative sequence
		constructed by the rule: \\\centerline{$x_{k+1} = x_{k} +
			\lambda_{k}s_{k}, \, k \in \mathbb{N},$}\smallskip \\ (f1)\enspace
		$\varepsilon_{0}
		> 0,$ $\varepsilon_{k+1} =
		\varepsilon_{k}\cdot(1-\beta)^{1-i_{k}},$ $k \in \mathbb{N}.$ \\
		Then it is fulfilled $\varepsilon_{k} \leq \bar{\varepsilon}, \,
		\forall k \in \mathbb{N}$,\, where $\bar{\varepsilon} = max \left
		\{\varepsilon_{0},\, \genfrac{}{}{}{0}{\mu}{1-\beta} \right \} >
		0.$
	\end{lemm}
	The next two universal theorems provide the possibility of
	evaluating the expected decrease of the objective function value when one uses
	the $\varepsilon_{k}-$normalized direction of descent and step-size
	chosen according to one of the above-described rules (Rule 1 or Rule 2). Here we speak of these next two theorems as universal ones in the sense that their formulation and proof do not depend on concrete algorithms. Indeed, there are important only the facts that all of the descent directions are normalized and the step-size is regulated by  Rule 1 or Rule 2.  For the proof details of these theorems,  we refer the interested reader to 
	\cite{Gabid_gab11}.
	\begin{theorem} \label{Gabid_t3.3} (estimate of the magnitude of
		decreasing the objective function value when the step length
		is selected according to Rule 1)(\cite{Gabid_gab11}, p.1089)\,\,
		If \\(b2)\enspace  the  conditions (b1), (c1), (e1) and (f1)  of Lemma
		\ref{Gabid_lem3.7}\, are fulfilled,\\
		(c2)\enspace \,$i_{k}$\, is the smallest index \,$i \in
		J(\hat{i})$\,
		for which there is fulfilled the condition (\ref{Gabid_e1})\, with $x = x_{k}$,\,
		$s = s_{k}$,\, $\varepsilon = \varepsilon_{k}$,\, $\lambda_{k} = \eta^{i_{k}}$,\,
		\,$\eta = (1 - \beta)^{1/(v - 1)}$,\, \,$\beta \in \left] 0,\,1\right[ $.\,
		
		Then there will be found a constant $\bar{C} > 0$\, such that for all
		\,$k \in \mathbb{N}$\, the following relation holds:
		\begin{multline} \label{Gabid_3.15}
		f(x_{k}) - f(x_{k + 1}) \geq -\bar{C} \cdot \langle
		\nabla f(x_{k}),s_{k}\rangle \geq  \\-\bar{C} \cdot (\langle
		\nabla f(x_{k}),s_{k}\rangle + \varepsilon_{k}
		\|s_{k}\|^{v}), 
	\end{multline}
			\begin{multline} \nonumber
	\mbox{where}\quad \,\bar{C} = \min \left\{C_{1},\,C_{2} \right\},\, \,C_{1} = \beta(1 - \beta)^{1/(v - 1)},\\ C_{2} = \left((1 - \beta)^{2}\varepsilon_{0}/\mu
		\right)^{1/(v - 1)}\beta. 
		\end{multline}
	\end{theorem}
	\begin{theorem} \label{Gabid_t3.4} (estimate of the magnitude of
		decreasing the objective function value when the step length
		is selected according to Rule 2)(\cite{Gabid_gab11}, p.1089)\,\,\,
		Let \\  (b3) \enspace the  conditions (b1), (c1), (e1) and (f1)  of Lem\-ma~\ref{Gabid_lem3.7}\, be fulfilled,\\
		(c3) \enspace the values of the iterative step-size \,$\lambda_{k}$,\,
		\,$\forall k \in \mathbb{N}$\, be determined using (\ref{Gabid_e2}).
		Then there will be found a constant $\bar{C} > 0$\, such that for all
		\,$k \in \mathbb{N}$\, the inequality (\ref{Gabid_3.15}) holds
		with the following constants:
		\begin{multline}\nonumber
		\bar{C} = \min \left\{C_{1},\,C_{2} \right\},\,C_{1} = \beta(1 - \beta)^{1/(v - 1)},\\
		\,C_{2} = \left(\varepsilon_{0}\mu^{-1} (1 - \beta)^{2}\right)^{1/(v - 1)}
		\left(1 + \varepsilon_{0}\mu^{-1}(1 - \beta)\beta^{-1}\right)^{-1}.\,
		\end{multline}
	\end{theorem}
		\begin{theorem} \label{Gabid_t3.11}(sublinear rate of convergence of
		ASDM)\,\,
		If \\ (b4)\enspace \,$f(x)$\, is a continuously differentiable
		pseudo-convex function on the convex set \,$D \subseteq
		\mathbb{R}^{n}$ ($X^{*}\subset D$)\, satisfying Condition A\, with some constant \,$\mu$\, and a~func\-tion \,$\tau(x,	y) = \|x - y\|^{v},$\, \,$v \geq 2,$\,   
		\\
		(c4)\enspace a numeric sequence
		$\{\theta_{k}\}$,\, which is defined by (\ref{Gabid_eq22}),\, satisfies the
		condition:   $\exists \theta
		> 0$\, such that \,$\theta_{k} \geq \theta$,\, $\forall k,$\\
		(d4)\enspace there exists a constant $\gamma > 0$ such that $\|\nabla f(x)\| \leq
		\gamma < \infty,$ $\forall x \in D,$ \\
		(e4)\enspace the Lebesgue set of the function \,$f(x)$\, at the point \,$x_{0}
		\in \mathbb{R}^{n},$\,  which is denoted by  \,$M(f,x_{0}) := \{x \in \mathbb{R}^{n}: f(x) \leq
		f(x_{0})\}$,\, is bounded,\\
		(g4)\enspace a step-size \,$\lambda_{k},\, k \in \mathbb{N}$\, is chosen
		according to one of the rules (Rule~1 or Rule~2).\\
		Then the sequence \,$\{x_{k}\}$,\, \,$k \in \mathbb{N}$\, 
		converges 	weakly , i.e. $$f(x_{k}) - f^{*} \sim O(1/k),$$ or
		equivalently, there exists a constant $C_{3} > 0$ such that it holds
		$$f(x_{k}) - f^{*} \leq C_{3} \cdot k^{-1}.$$
	\end{theorem}
	\proof
		Clearly, \,$f(p^{*}_{k}) = f^{*},$\, \,$f(x_{k}) > f(p^{*}_{k}),$\,
		\,$\forall k \in \mathbb{N}.$\, By definition of pseudo-convex functions, it
		is fulfilled \,$\langle \nabla f(x_{k}), p^{*}_{k} - x_{k} \rangle < 0.$\,
		According to the assertions of Theo\-rems~\ref{Gabid_t3.3}--\ref{Gabid_t3.4},\,
		regardless
		of whether Rule 1 will be selected for calculating
		the step-size or Rule 2, there may be found
		a constant \,$\bar{C} > 0$\, such that 
		the inequality (\ref{Gabid_3.15})\, holds for all \,$k \in \mathbb{N}$.\, Choose the
		subset of indices \,$\mathbb{N}_{1} \subset \mathbb{N}$\, such that
		\,$s_{k} = p_{k},$\, \,$k \in \mathbb{N}_{1}.$\,   One  then has
		\,$s_{k} = \genfrac{}{}{}{0}{p_{k}}{\varepsilon_{k}\|p_{k}\|^{v-2}}$\,  \,$\forall k
		\in \mathbb{N}_{2} = \mathbb{N} \backslash \mathbb{N}_{1}$.\, For all $k \in
		\mathbb{N}_{1}$, we have the estimate:
		\begin{equation}\nonumber
		f(x_{k}) - f(x_{k+1}) \geq -\bar{C}\cdot \langle \nabla f(x_{k}), s_{k} \rangle = \bar{C}\cdot \|\nabla f(x_{k})\|^{2}.
		\end{equation}
		From Lemma \ref{Gabid_lem3.7}, due to (\ref{Gabid_3.15}), for all
		$k \in \mathbb{N}_{2}$,\, it follows the relation
		\begin{multline}\nonumber
		f(x_{k}) - f(x_{k+1}) \geq  -\bar{C}\cdot \langle \nabla f(x_{k}), s_{k} \rangle =\\
		\genfrac{}{}{}{0}{-\bar{C}}{\varepsilon_{k}\|p_{k}\|^{v-2}}\cdot \langle \nabla f(x_{k}), p_{k} \rangle  \geq
		\genfrac{}{}{}{0}{\bar{C}}{\bar{\varepsilon}
			\gamma^{v-2}}\|\nabla f(x_{k})\|^{2}.
		\end{multline}
		Thus, for all $k \in \mathbb{N}$, one has arrived at the
		inequality
		\begin{equation}\label{Gabid_3.30}
		f(x_{k}) - f(x_{k+1}) \geq
		\tilde{C}\cdot
		\|\nabla f(x_{k})\|^{2},
		\end{equation}
		where \,$\tilde{C} =
		\bar{C} \min
		\left\{1 ,\,
		\genfrac{}{}{}{0}{1}{\bar{\varepsilon}\gamma^{v-2}}\right \}.$\,
		By virtue of the condition (e4), $$diam\, M(f,x_{0}) = \sup \{ \|x - y\|, \,x,\,y \in M(f,x_{0}) \} = \bar{\eta} < + \infty.$$
		Therefore, we immediately have the following estimate 
		\begin{equation}\nonumber
		\|\nabla f(x_{k})\|^{2} \geq \|\nabla f(x_{k})\|^{2} \|x_{k} - p^{*}_{k}\|^{2} \cdot \dfrac{1}{\bar{\eta}^{2}} \geq 
		\dfrac{1}{\bar{\eta}^{2}} \cdot \langle \nabla f(x_{k}), x_{k} - p^{*}_{k} \rangle^{2}.
		\end{equation}
		Consequently, for all \,$k \in \mathbb{N}$\, from (\ref{Gabid_eq22}),\,(\ref{Gabid_3.30}),\, and (c4) it follows
		\begin{equation}\label{Gabid_3.337}
		f(x_{k}) - f(x_{k+1}) \geq
		\dfrac{\tilde{C}}{\bar{\eta}^{2}}\cdot
		\langle \nabla f(x_{k}), x_{k} - p^{*}_{k} \rangle^{2} = C_{3} \cdot ( f(x_{k}) - f(p^{*}_{k}))^{2},
		\end{equation}
		where \,$C_{3} = \theta^{2}\dfrac{\tilde{C}}{\bar{\eta}^{2}}.$\,
		Due to Lemma \ref{Gabid_lem2.5},\, the latter implies that the
		sequence $\{x_{k}\},\, k \in \mathbb{N}$ is weakly convergent to
		a solution of (\ref{Gabid_eq1}),\, since there holds the following
		estimate for the convergence rate:
		$$f(x_{k}) - f^{*} \leq
		C_{3}^{-1}k^{-1}. \quad \quad \qquad \enspace \square$$
		 
			\begin{remark}
		There is no a need to strictly require the fulfillment of the condition (b4) of Theorem \ref{Gabid_t3.11}\, in the whole space \,$\mathbb{R}^{n}.$\, For instance, it is sufficient  to have  that 
		\,$D = M(f, x_{0}),$\, where \,$x_{0}$\, is a starting point for ASDM. Besides, there may simply be chosen some convex set for which it holds  \,$X^{*} \subset D$.\, 
	\end{remark}
	Preliminary computational tests confirm the efficiency of the proposed method and the strict monotone property of the used step-size rules. 
	These tests show that the results of minimization depends on the user-selected parameters of ASDM such as \,$\beta$\, and \,$\varepsilon$.\, We also observe that at each iteration (after the first one) there is needed 
	only one function and gradient evaluation. Our preliminary experiments
	has demonstrated the ability  of ASDM to lead to the minimum neighborhood at low computational costs.
	
	\section{Conclusions} \label{Gabid_sec6}
	We proposed a fully adaptive variant of the steepest descent method. There are used some novel rules for the calculation of the step length in which the iteration step-size is regulated additionally by an
	adaptation of the $\varepsilon-$normalization parameter for the
	descent direction. The finiteness of all the
	procedures of adaptive cont\-rol\-ling both the
	parameter of an~$\varepsilon-$nor\-ma\-lization of a descent
	direction and a step length was established. For the problem of unconstrained minimizing a smooth pseudo-convex function, we justified the sublinear rate of the convergence for the adaptive variant of the steepest descent method.
	
	One of the motivating ideas was that of using in the future the adaptive steepest descent method to solve the problems of sets separation (by minimizing the error function) and related problems of data mining (in particular, neural network classification of data).
	
	%\begin{acknowledgements}
	%	The author thanks the anonymous referees and the editor for their helpful comments %\end{acknowledgements}
	% Insert your figure, if needed.
	%\begin{figure}[!h]
	%  \centering
	% \includegraphics[width=0.7\maxpicturewidth]{your_figure.eps}
	%  \center{Fig.~1. Your caption.}
	%\end{figure}
	
\end{document}